\documentclass[11pt]{amsart}
\usepackage{amsmath,amssymb, graphics, amscd,latexsym}

\makeatletter

\newtheorem{Theorem}{Theorem}
\newtheorem{Lemma}[Theorem]{Lemma}
\newtheorem{Corollary}[Theorem]{Corollary}
\newtheorem{Proposition}[Theorem]{Proposition}
\newtheorem{Claim}[Theorem]{Claim}
\newtheorem{Definition}[Theorem]{Definition}

\theoremstyle{remark}
\newtheorem{Remark}[Theorem]{Remark}

\newtheorem{Conjecture}[Theorem]{Conjecture}

\begin{document}
\newcommand{\eps}{\varepsilon}
\newcommand{\om}{\omega}
\newcommand\Om{\Omega}
\newcommand\la{\lambda}
\newcommand\vphi{\varphi}
\newcommand\vrho{\varrho}
\newcommand\al{\alpha}
\newcommand\La{\Lambda}
\newcommand\si{\sigma}
\newcommand\be{\beta}
\newcommand\Si{\Sigma}
\newcommand\ga{\gamma}
\newcommand\Ga{\Gamma}
\newcommand\de{\delta}
\newcommand\De{\Delta}

\newcommand\cA{\mathcal  A}
\newcommand\cB{\mathcal B}
\newcommand\cD{\mathcal  D}
\newcommand\cM{\mathcal  M}
\newcommand\cN{\mathcal  N}
\newcommand\cT{\mathcal  T}
\newcommand\cP{\mathcal  P}
\newcommand\cp{\mathcal  p}
\newcommand\cQ{\mathcal  Q}
\newcommand\cG{\mathcal G}
\newcommand\cq{\mathcal  q}
\newcommand\cc{\mathcal  c}
\newcommand\cs{\mathcal  s}
\newcommand\cS{\mathcal  S}
\newcommand\ct{\mathcal  t}
\newcommand\cZ{\mathcal  Z}
\newcommand\cR{\mathcal  R}
\newcommand\cu{\mathcal  u}
\newcommand\cU{\mathcal  U}
\newcommand\cI{\mathcal  I}
\newcommand\cJ{\mathcal  J}
\newcommand\co{\mathcal  o}
\newcommand\cO{\mathcal  O}
\newcommand\cv{\mathcal  v}
\newcommand\cV{\mathcal  V}
\newcommand\cx{\mathcal  x}
\newcommand\cX{\mathcal  X}
\newcommand\cw{\mathcal  w}
\newcommand\ck{\mathcal  k}
\newcommand\cK{\mathcal  K}
\newcommand\cW{\mathcal  W}
\newcommand\cz{\mathcal  z}
\newcommand\cy{\mathcal  y}
\newcommand\ca{\mathcal  a}
\newcommand\ch{\mathcal  h}
\newcommand\cH{\mathcal  H}
\newcommand\cF{\mathcal F}
\newcommand\bfG{\mbox {\bf  G}}
\newcommand\bfg{\mbox {\bf  g}}
\newcommand\bfC{\mbox {\bf  C}}
\newcommand\bfN{\mbox {\bf  N}}
\newcommand\bfT{\mbox {\bf  T}}
\newcommand\bfP{\mbox {\bf  P}}
\newcommand\bfp{\mbox {\bf  p}}
\newcommand\bfQ{\mbox {\bf  Q}}
\newcommand\bfq{\mbox {\bf  q}}
\newcommand\bfc{\mbox {\bf  c}}
\newcommand\bfs{\mbox {\bf  s}}
\newcommand\bfS{\mbox {\bf  S}}
\newcommand\bft{\mbox {\bf  t}}
\newcommand\bfZ{\mbox {\bf  Z}}
\newcommand\bfR{\mbox {\bf  R}}
\newcommand\bfu{\mbox {\bf  u}}
\newcommand\bfU{\mbox {\bf  U}}
\newcommand\bfo{\mbox {\bf  o}}
\newcommand\bfO{\mbox {\bf  O}}
\newcommand\bfv{\mbox {\bf  v}}
\newcommand\bfV{\mbox {\bf  V}}
\newcommand\bfx{\mbox {\bf  x}}
\newcommand\bfX{\mbox {\bf  X}}
\newcommand\bfw{\mbox {\bf  w}}
\newcommand\bfk{\mbox {\bf  k}}
\newcommand\bfK{\mbox {\bf  K}}
\newcommand\bfW{\mbox {\bf  W}}
\newcommand\bfz{\mbox {\bf  z}}
\newcommand\bfy{\mbox {\bf  y}}
\newcommand\bfa{\mbox {\bf  a}}
\newcommand\bfh{\mbox {\bf  h}}
\newcommand\bfH{\mbox {\bf  H}}
\newcommand\bfJ{\mbox {\bf  J}}
\newcommand\bfj{\mbox {\bf  j}}
\newcommand\bbC{\mbox {\mathbb C}}
\newcommand\bbN{\mbox {\mathbb N}}
\newcommand\bbT{\mbox {\mathbb T}}
\newcommand\bbP{\mbox {\mathbb P}}
\newcommand\bbQ{\mbox {\mathbb Q}}
\newcommand\bbS{\mbox {\mathbb S}}
\newcommand\bbZ{\mbox {\mathbb Z}}
\newcommand\bbR{\mbox {\mathbb R}}
\newcommand\bbU{\mbox {\mathbb U}}
\newcommand\bbO{\mbox {\mathbb O}}
\newcommand\bbV{\mbox {\mathbb V}}
\newcommand\bbX{\mbox {\mathbb X}}
\newcommand\bbK{\mbox {\mathbb K}}
\newcommand\bbW{\mbox {\mathbb W}}
\newcommand\bbH{\mbox {\mathbb H}}

\newcommand\apeq{\fallingdotseq}
\newcommand\Lrarrow{\Leftrightarrow}
\newcommand\bij{\leftrightarrow}
\newcommand\Rarrow{\Rightarrow}
\newcommand\Larrow{\Leftarrow}
\newcommand\nin{\noindent}
\newcommand\ninpar{\par \noindent}
\newcommand\nlind{\nl \indent}
\newcommand\nl{\newline}
\newcommand\what{\widehat}
\newcommand\tl{\tilde}
\newcommand\wtl{\widetilde}
\newcommand\order{\mbox{\text{order}\/}}
\newcommand\GL{\text{GL}\/}
\newcommand\PGL{\text{PGL}\/}
\newcommand\Spec{\text{Spec}\/}
\newcommand\weight{\text{weight}\/}
\newcommand\ord{\text{ord}\/}
\newcommand\Int{\text{Int}\/}
\newcommand\grad{\text{grad}\/}
\newcommand\Ind{\text{Ind}\/}
\newcommand\Disc{\text{Disc}\/}
\newcommand\Ker{\text{Ker}\/}
\newcommand\Image{\text{Image}\/}
\newcommand\Coker{\text{Coker}\/}
\newcommand\Id{\text{Id}\/}
\newcommand\id{\text{id}}
\newcommand\dsum{\text{\amalg}}
\newcommand\val{\text{val}}
\newcommand\mv{\text{m-vector}}
\newcommand\iv{\text{i-vector}}
\newcommand\minimum{\text{minimum}\/}
\newcommand\modulo{\text{modulo}\/}
\newcommand\Aut{\text{Aut}\/}
\newcommand\PSL{\text{PSL}\/}
\newcommand\Res{\text{Res}\/}
\newcommand\rank{\text{rank}\/}
\newcommand\codim{\text{codim}\/}
\newcommand\msdim{\text{ms-dim}\/}
\newcommand\icodim{\text{i-codim}\/}
\newcommand\ocodim{\text{o-codim}\/}
\newcommand\emdim{\text{ems-dim}\/}

\newcommand\Cone{\text{Cone}\/}
\newcommand\maximum{\text{maximum}\/}
\newcommand\Vol{\text{Vol}\/}
\newcommand\Coeff{\text{Coeff}\/}
\newcommand\lcm{\text{lcm}\/}
\newcommand\degree{\text{degree}\/}
\newcommand\Pol{\cal {POL}}
\newcommand\Ts{Tschirnhausen}
\newcommand\TS{Tschirnhausen approximate}
\newcommand\Stab{\text{Stab}\/}
\newcommand\civ{complete intersection variety}
\newcommand\nipar{\par \noindent}
\newcommand\wsim{\overset{w}{\sim}}
\newcommand\Gr{\bfZ_2*\bfZ_3}
\newcommand\QED{~~Q.E.D.}
\newcommand\bsq{$\blacksquare$}
\newcommand\bff{\mbox {\bf  f}}
\newcommand\newcommandby{:=}
\newcommand\inv{^{-1}}
\newcommand\nnt{(\text{nn-terms})}
\renewcommand{\subjclassname}{\textup{2000} Mathematics Subject Classification}
\def\mapright#1{\smash{\mathop{\longrightarrow}\limits^{{#1}}}}
\def\mapleft#1{\smash{\mathop{\longleftarrow}\limits^{{#1}}}}
\def\mapdown#1{\Big\downarrow\rlap{$\vcenter{\hbox{$#1$}}$}}
\def\mapdownn#1#2{\llap{$\vcenter{\hbox{$#1$}}$}\Big\downarrow\rlap{$\vcenter{\hbox{$#2$}}$}}
\def\mapup#1{\Big\uparrow\rlap{$\vcenter{\hbox{$#1$}}$}}
\def\mapupp#1#2{\llap{$\vcenter{\hbox{$#1$}}$}\Big\uparrow\rlap{$\vcenter{\hbox{$#2$}}$}}

\def\rdown#1{\searrow\rlap{$\vcenter{\hbox{$#1$}}$}}
\def\semap#1{\searrow\rlap{$\vcenter{\hbox{$#1$}}$}}

\def\rup#1{\nearrow\rlap{$\vcenter{\hbox{$#1$}}$}}
\def\nemap#1{\nearrow\rlap{$\vcenter{\hbox{$#1$}}$}}
\def\ldown#1{\swarrow\rlap{$\vcenter{\hbox{$#1$}}$}}
\def\swmap#1{\swarrow\rlap{$\vcenter{\hbox{$#1$}}$}}
\def\lup#1{\nwarrow\rlap{$\vcenter{\hbox{$#1$}}$}}
\def\nwmap#1{\nwarrow\rlap{$\vcenter{\hbox{$#1$}}$}}
\def\defby{:=}
\def\eqby#1{\overset {#1}\to =}
\def\inv{^{-1}}
\def\bnu{{(\nu)}}
\def\ocup#1{\underset{#1}\cup}
\title[ Alexander polynomial of sextics 
]
{Alexander polynomial of sextics }

\author
[M. Oka ]
{Mutsuo Oka }
\address{\vtop{
\hbox{Department of Mathematics}
\hbox{Tokyo Metropolitan University}
\hbox{1-1 Mimami-Ohsawa, Hachioji-shi}
\hbox{Tokyo 192-0397}
\hbox{\rm{E-mail}: {\rm oka@comp.metro-u.ac.jp}}
}}

\begin{abstract}
 Alexander polynomials  of sextics are computed in the case
of sextics  with  only simple singularities
or sextics of torus type with arbitrary singularities.
We will show that for irreducible sextics,
there are 
only 4 possible Alexander polynomials:
 $(t^2-t+1)^j,\,j=0,1,2,3$.
For the computation, we use the method of Esnault-Artal
and the classification result in  our previous papers.
\end{abstract}
\maketitle

\pagestyle{headings}
\section{Introduction}
Recall that a sextics is called  {\em of torus type} if it is defined by 
polynomial 
$f(x,y):=f_2(x,y)^3+f_3(x,y)^2$ where $\degree f_i(x,y)=i$ for $i=2,3$.
Sextics without such a torus expression is called {\em of non-torus type}.
In \cite{Oka-Pho1}, we have computed the fundamental group
of the complement of tame sextics of torus type
and we have classified the  configurations of the singularities
on non-tame sextics of torus type in \cite{Oka-Pho2}.
In general, the computation of the
fundamental groups is  more complicated in the case of
non-tame torus curves.
In this paper, we are interested in another invariant which is  called
{\em an  Alexander polynomial}. This invariant is weaker than the
fundamental groups,
but easier to be computed.
In fact,  there are Zariski pairs with different
fundamental group but with the same Alexander polynomial (\cite{Two}).
The advantage of  Alexander polynomial is that it can be computed by
the  data of the local singularities and the data about their  global position
in $\bfP^2$.

It is the purpose of this paper to give a complete atlas of the
Alexander polynomials of \nl
(a) sextics of torus type with arbitrary
singularities
or \nl
 (b) sextics
 with only simple singularities, not necessarily of torus type.

In fact, we will show
that 
 there are only four possibilities of Alexander polynomials:
\[
 1,\quad (t^2-t+1),\quad (t^2-t+1)^2,\quad (t^2-t+1)^3.
\]
 for irreducible sextics of type (a) or (b).
The case $\De(t)=(t^2-t+1)^3$ corresponds to 9 cuspidal sextics,
$\De(t)=1$ corresponds to sextics of non-torus type with $\rho(5)\le 6$.
It is expected that every sextics of non-torus type satisfies
the above inequality. The case
$\De(t)=(t^2-t+1)^2$ corresponds to the cases:
\[
 \Si(C)=[8A_2],\, [8A_2,A_1],\, [6A_2,E_6]
,\, [6A_2,A_5],\,
[6A_2,A_5,A_1],\, [4A_2,2A_5],\,
[4A_2,A_5,E_6] \]
and  sextics of torus type with  configuration
    $ [B_{3,6},3A_2],\, [C_{3,9},3A_2]$.
For reducible sextics of torus type, we can have further
\[
\tilde\De(t)=(t^2-t+1)^2(t^2+t+1),\,
(t^2-t+1)(t^2+t+1),\,\]
\[
 (t^2-t+1)^2(t^2+t+1)(t+1)^2,\,
(t^2-t+1)^4(t^2+t+1)^4(t+1)^4\, 
\]
The last case corresponds to a  sextic with
$B_{6,6}$ and $C$ has 6 line components.
See Corollary \ref{Simple-list} and Corollary \ref{NonSimple-list}
for further detail.

Sextics of non-torus type with non-simple singularities are not
considered in this paper.
\section{Alexander Polynomial of a plane curve}
Let $C$ be an affine curve defined by a polynomial $f(x,y)$ of degree $d$.
We assume that the line at infinity $L_\infty$ and $C$ intersect
transversely and we identify $\bfP^2-L_\infty\cong \bfC^2$.
Let $\phi: \pi_1(\bfC^2-C) \to \bfZ$ be the canonical homomorphism
induced by the composition of the Hurewicz homomorphism 
$\psi:\pi_1(\bfC^2-C)\to H_1(\bfC^2-C,\bfZ)\cong  \bfZ^r$
and the
summation
homomorphism
$\bfZ^r\to \bfZ$  where 
$r$ is the number of irreducible components of $C$.
Let $t$ be a generator of $\bfZ$.

Let $\tilde X\to \bfC^2-C$ be the infinite cyclic covering corresponding 
to $\Ker \phi$. Then $H_1(\tilde X,\bfQ)$ has the structure of $\La$
module
where
$\La:=\bfQ[t,t\inv]$. Thus we can write as 
\[
 H_1(\tilde X,\bfQ)\cong \bfQ^{r-1}\oplus \La/\la_1(t)\oplus\cdots\oplus \La/\la_\nu(t)
\]
where  $\la_i(t)$ is a polynomial with integral
coefficients satisfying the properties
 $\la_i(0)\ne 0$ and  $\la_i|\la_{i+1}$  for $i=1,\dots, \nu-1$. The 
{\em generic Alexander polynomial} $\De(t)$ is defined by the product
$\De(t)=(t-1)^{r-1}\la_1(t)\cdots \la_\nu(t)$.
Let us call $\tilde\De(t):=\De(t)/(t-1)^{r-1}$
{\em the reduced Alexander polynomial}. It does not depend
on the choice of the generic line at infinity $L_\infty$. Let $F(X,Y,Z)$ be the defining homogeneous polynomial of $C$ and 
let $M=F\inv(1)\subset \bfC^3$ be the Milnor fiber of $F$.
Randell showed in \cite{Randell}
 that $\De(t)$ can be computed as the characteristic 
polynomial of the monodromy on th first homology group
 of the Milnor  fibration defined by 
the homogeneous polynomial $F(X,Y,Z)$. Thus the degree of  $\De(t)$ is
equal to the first Betti number $b_1(M)$.
On the other hand, Libgober has proved that the degree of  
$\tilde\De(t)$ is  equal to
the sum $\sum_{j=1}^{d-1}\beta_j$, where $\beta_j$ is defined by the number of 
factors  $\ell,\ell=1,\dots, \nu$, such that $\exp(2\pi j \sqrt{-1}/d)$ is 
a root of $\la_\ell(t)=0$.
Combining these results, we observe that {\em $\la_j(t)$ has no multiple roots}. 
Furthermore Libgober and Loeser-Vaqui\'e  showed that 
\begin{Lemma} {\rm (\cite{LibgoberArcata},\cite{Loeser-Vaquie})}
The polynomial 
$\De(t)$ is  written as the product 
\[
\tilde \De(t)=\prod_{k=1}^{d-1}\De_k(t)^{\ell_k},\quad k=1,\dots,d-1
\]
where 
\[
\De_k=(t-\exp(\frac {2 k \pi i}d))(t-\exp(\frac {-2 k \pi i}d))
\]
and  $\ell_k=\dim H^1(P^2,{L^{(k)}})$.
\end{Lemma}
For the definition of the sheaf $L^{(k)}$,
we refer to \cite{Esnault, Artal,Loeser-Vaquie}.
We use the method of Esnault-Artal 
to compute $\ell_k$ (\cite{Esnault}, \cite{Artal}).
 Note that for the case of sextics $d=6$, 
$\De_5(t)=\De_1(t)=t^2-t+1$, $\De_4(t)=\De_2(t)=t^2+t+1$
and $\De_3(t)=(t+1)^2$.

Let $C$ be a given plane curve of degree $d$ defined by $f(x,y)=0$
 and let $\Si( C)$ be the singular locus of $C$ and let $P\in \Si (C)$
be a singular point. Consider  an embedded resolution of $C$, $\pi:\tilde U\to U$ 
 and let $E_1,\dots, E_s$ be the exceptional divisors. 
Let us choose $(u,v)$ be a local coordinate system centered at $P$
and let $k_i$ and $m_i$ be the order of zero of 
the canonical two form $\pi^*(du\land dv)$ 
and $\pi^*f$ respectively along  the divisor $E_i$. An ideal
 $\cJ_{P,k,d}$ of $\cO_{P}$
is generated by the function germ $\phi$ such that 
the pull-back $\pi^*\phi$  vanishes of order at least $-k_i+[km_i/d]$
along $E_i$. Namely
\[
 \cJ_{P,k,d}=\{\phi\in \cO_{P}; (\pi^*\phi)\ge \sum_{i}(-k_i+[km_i/d])E_i\}
\]
Let us consider
 the canonical homomorphisms induced by the restrictions:
\begin{eqnarray*}
 \si_{k,P}: \cO_{P}\to \cO_{P}/\cJ_{P,k,d},\quad
 \si_{k}: H^0(\bfP^2,\cO(k-3))\to \bigoplus_{P\in \Si (C)}\cO_{P}/\cJ_{P,k,d}
\end{eqnarray*}
where the right side of $\si_k$ is the sum over singular points of $C$.
By \cite{Artal},
the integer $\ell_k$ is equal to the dimension of the 
cokernel $\si_{k}$.
\subsection{Non-degenerate singularity} 
In general, the computation of the kernel and cokernel of the homomorphism 
$\si_{k,P}$ requires   an explicit computation of the resolution.  
However for  the class of non-degenerate singularities, these informations 
can be obtained easily by a toric resolution.
See \cite{Var-zeta,Okabook} for the definition of non-degenerate 
singularities. In fact, we 
do not need  any data  in detail for a resolution.
Let us assume that $(u,v)$ is a fixed  local analytic coordinate such
that 
$P=(0,0)$ and $C$ is defined by a  
function germ $f(u,v)$ and the Newton boundary $\Ga(f)$ is non-degenerate with respect to $(u,v)$. 
Let $Q_1,\dots,Q_s$ be the primitive covectors which correspond to the
edges
$\De_1,\dots,\De_s$ of $\Ga(f)$
respectively. That is, $\De(Q_i,f)=\De_i$.
Here we use the same terminology as in \cite{Okabook}. 

 Let $\pi:\tilde U\to U$ is the canonical toric modification 
and let $\hat E(Q_i)$ be the exceptional divisor corresponding to $Q_i$. Recall that the order of zeros 
of the canonical two form $\pi^*(du\land dv)$ along
 the divisor $\hat E(Q_i)$ is simply given by  
$|Q_i|-1$ where $|Q_i|:=(p+q)$ for $Q_i={}^t(p,q)$.  See \cite{Okabook},
p. 178.
 For a function germ $g(u,v)$, 
let $m(g,Q_i)$ be the multiplicity of the  
pull-back $\pi^*g$ on $\hat E(Q_i)$.  This number is equal to
 $d(Q_i,g)$ in \cite{Okabook}, which is the minimal value of the linear
 function $Q_i$ restricted to the Newton boundary of $g(u,v)$.
Then the following criterion is essentially due to Merle-Teissier \cite{Merle-Teissier}  and it is useful for the computation of  
the ideal $\cJ_{P,k,d}$. 

\begin{Lemma}\label{non-degenerate-rho} Take a function germ $g$. 
Then  $g$ is contained
  in the ideal $\cJ_{P,k,d}$ 
if and only if  
\[(\sharp)\quad m(g,Q_i)\ge [\frac kd m(f,Q_i)]-|Q_i|+1,\quad\text{for}\,\, 
i=1,\dots,s\] 
The ideal $\cJ_{P,k,d}$ is generated by the monomials satisfying the above condition. 
\end{Lemma} 
Note that the condition $(\sharp)$ can be checked without constructing
an explicit toric compactification.

{\em Proof.} Let $\Si^*$ be a regular subdivision of $\Ga^*(f)$ which we
use
to construct our toric modification $\pi: \tilde U \to U$
We put  $Q_0={}^t(1,0)$ and $Q_{s+1}={}^t(0,1)$. 
Then for $Q_0,Q_{s+1}$ the conditions are satisfied. Take the cone 
$\text{Cone}({Q_i,Q_{i+1}})$ 
and let $T_{i,1},\dots, T_{i,\nu_i}$ be the covectors which are 
inserted for the regular subdivision 
of the cone $\Cone(Q_i,Q_{i+1})$.  
Note that the conditions $(\sharp)$ for $Q_{i}, Q_{i+1}$ are equivalent to 
\begin{eqnarray*} 
m(Q_j,g)> \frac kd m(Q_j,f)-|Q_j|,\quad j=i,i+1 
\end{eqnarray*} 
On the other hand, we can write $T_{i,j}=\al_j Q_i+\be_j Q_{i+1}$ for some positive rational numbers $\al_j,\be_j$. 
We also observe that 
\begin{eqnarray*} 
m(T_{i,j},f)=\al_j m(Q_i,f)+\be_j m(Q_{i+1},f)\\ 
m(T_{i,j},g)\ge \al_j m(Q_i,g)+\be_j m(Q_{i+1},g). 
\end{eqnarray*} 
The first equality follows from the property:
$\De(T_{ij},f)=\De_i\cap \De_{i+1}$.
As $|T_{i,j}| =\al_j|Q_i|+\be_i |Q_{i+1}|$, we obtain 
\begin{multline*} 
\lefteqn{ m(T_{i,j},g)\ge \al_j m(Q_i,g)+\be_j m(Q_{i+1},g)} \\ 
>\frac kd (\al_j m(Q_i,f)+\be_j m(Q_{i+1},f))-(\al_j|Q_i|+\be_j|Q_{i+1}|)
=\frac kd m(T_{i,j},f)-|T_{i,j}| 
\end{multline*}
Thus the condition $(\sharp)$ is satisfied for the  covectors $T_{ij}$. 
The second assertion is obvious. \qed

\section{Alexander polynomials of Sextics with simple singularities} 
 \subsection{Normal forms of simple singularities}
For simple singularities, we use the following normal forms.
\begin{eqnarray}\label{normal-form1}
 \begin{cases} 
&A_n:  f(u,v)=v^2+u^{n+1}+\text{(higher terms)}\\
& D_k: f(u,v)=v^2 u+u^{k-1}+\text{(higher terms)}\\
&E_6: f(u,v)=v^3+u^4+\text{(higher terms)}\\
&E_7: f(u,v)=v^3+vu^3+\text{(higher terms)} \\ 
&E_8: f(u,v)=v^3+u^5+\text{(higher terms)} 
\end{cases}
\end{eqnarray}
Note that $\deg f(0,v)\ge 3$  for $D_j$ and $\deg f(u,0)\ge 5$ for
$E_7$ and in
 fact we can make $f(x,y)$ to be convenient  by
 taking $u\mapsto u+ v$ or $v\mapsto v+u^2$.
In the above notation, we observe that the line $v=0$
gives the tangent cone of the respective singularities and
 the local intersection number with the given curve 
 is strictly greater than
the respective multiplicity. In the case of $D_k$ singularity, $v=0$ is 
one of the tangent cone. Let $L_v$ be the projective line passing
through the origin $O=(0,0)$ in the coordinate $(u,v)$ and tangent to
the
smooth curve $v=0$.  We call $L_v$ {\em the principal tangential direction}
of
the simple singularity. As $L_v$ is written as
$v+a\, u^2+\text{(higher terms)}$,
we have  $I(L_v,C;O)\ge 4$ for $A_j,\, j\ge 3$ or $E_6,E_7,\, E_8$ or
$D_k,\, k\ge 5$.
Here $I(f,g;P)$ is the local intersection number of $f=0$ and $g=0$ at $P$.
Let us consider the canonical homomorphisms defined in the previous section:
\begin{eqnarray*}
\si_{k,P}: \,\cO_{P}\to \cO_{P}/\cJ_{P,k,6},\qquad
\si_{k}&: \,H^0(\bfP^2,\cO(k-3))\to
\bigoplus_{P}\cO_{P}/\cJ_{P,k,6}
\end{eqnarray*}

\subsection{Description of local data for simple singularities}  
We assume hereafter $(u,v)$ is a chosen local coordinate system so that
the defining equations are written as in (\ref{normal-form1}).
The local data for the simple singularities   
are described by the following. 
\begin{Proposition}\label{simple-rho}  
Assume that $(C,P)$ is a simple singularity defined by a normal form as
 in (\ref{normal-form1}). 
For $k\le 3$, $\rho(P,k)=0$.
For $k=4$ we have  
\begin{eqnarray*}  
\rho(P,4)&=&0,~ \text{for}~(C,P)\cong A_1,\dots,A_4 \\  
\rho(P,4)&=&1,~ \cJ_{P,4,6}=  
\langle u,v\rangle\quad \text{for}~(C,P)\cong A_5,\dots,A_{10}, D_4,\dots,D_9,E_6, \,E_7, \,E_8 \\  
\rho(P,4)&=&2,~ \cJ_{P,4,6}=  
\langle u^2,v\rangle\quad \text{for}~(C,P)\cong A_{11},\dots,A_{16}, D_{10},\dots,D_{15} \\  
\rho(P,4)&=&3,~ \cJ_{P,4,6}=  
\langle u^3,v\rangle\quad \text{for}~(C,P)\cong A_{17},\dots,A_{22}, D_{16},\dots,D_{21}  
\end{eqnarray*}  
For $k=5$, we have $\rho(P,5)=0$ for $A_1$ and   
\begin{eqnarray*}  
\rho(P,5)&=&1,~  
\cJ_{P,5,6}=  
\langle u,v\rangle~\text{for}~(C,P)\cong A_2,A_3,A_4,D_4,D_5 \\  
\rho(P,5)&=&2,~  
\cJ_{P,5,6}=  
\langle u^2,v\rangle~\text{for}~(C,P)\cong A_5,A_6,A_7,E_6,E_7,E_8,D_6,D_7,D_8 \\  
\rho(P,5)&=&3,~  
\cJ_{P,5,6}=  
\langle u^3,v\rangle~\text{for}~(C,P)\cong A_8,A_9,A_{10},D_9,D_{10},D_{11} \\  
\rho(P,5)&=&4,~  
\cJ_{P,5,6}=  
\langle u^4,v\rangle~\text{for}~(C,P)\cong A_{11},A_{12},A_{13},
D_{12},D_{13},D_{14} \\  
\rho(P,5)&=&5,~  
\cJ_{P,5,6}=  
\langle u^5,v\rangle~\text{for}~(C,P)\cong A_{14},A_{15}, A_{16},D_{15},
D_{16},D_{17} \\  
\rho(P,5)&=&6,~  
\cJ_{P,5,6}=  
\langle u^6,v\rangle~\text{for}~(C,P)\cong A_{17},A_{18},A_{19},D_{18},D_{19},
D_{20}  
\end{eqnarray*}  
\end{Proposition}  
We call  a singularity $(C,P)$ is
 {\em $\rho(k)$-essential} if $\rho(k,P)>0$.
Thus we see that
$A_1,\dots,A_4$  are not $\rho(4)$-essential
and  $A_1$ is not  $\rho(5)$-essential. 

\begin{Corollary}\label{Local-rho} 
{\rm I.} Take
a germ $g(x,y)\in \cO_P$.
Then $g$ is contained in the kernel of $\si_{k,P}$
if and only if
$I(g,v;P)\ge \rho(P,k)$.
 This implies also 
\begin{eqnarray}\label{int-mult}
 I(g,f;P)\ge 2\rho(P,k)
\end{eqnarray}

For $k=4$, the equality in  (\ref{int-mult}) 
holds only if $(C,P)\cong A_j,\, j\ge 5$.

For $k=5$,  the equality  in (\ref{int-mult})  holds only if
either $(C,P)\cong A_j,~j\ge 2$ or $(C,P)\cong E_6$
 singularities.
If  $(C,P)\cong E_6$ and  $I(g,f;P)=4$, the curve
 $g=0$ is smooth at $P$
and it is tangent to the tangent cone of $(C,P)$.

{\rm II.} Assume that $g,g'\in \Ker\, \si_{k,P}$. Then
$I(g,g';P)\ge \rho(P,k)$.
\end{Corollary}
{\em Proof.} Recall that  $f=0$ is the defining polynomial of $C$.
Let $(u,v)$ be a   coordinate system which gives the normal form 
(\ref{normal-form1}).
Write  $g=v\, g_1+u^r\, g_2$, $g_1,g_2\in \cO_P$ with
$g_2(0)\ne 0$.
Then  $g$ is contained in the ideal
$\cJ_{P,k,6}$ if and only if $r\ge \rho(P,k)$.
Using the normal form in (\ref{normal-form1}), we have
\begin{eqnarray}
I(u,f;P)& =&\begin{cases} 2\quad &(C,P)\cong A_j,j \ge 2\\
             3\quad &(C,P)\cong E_6,E_7,E_8               \\
\end{cases}
\label{int-u}\\
I(v,f;P)&=&\begin{cases} j+1\quad &(C,P)\cong A_j,j \ge 2\\
 4,\quad &(C,P)\cong E_6\\
\ge 5 \quad &(C,P)\cong E_7,\,E_8\\
             j-1\quad &(C,P)\cong D_j,j\ge 4
\end{cases}\label{int-v}
\end{eqnarray}
Thus if $r\ge \rho(P,k)$, we have
\begin{eqnarray*}
&I(u^r,f;P)\ge 2\, \rho(P,k),\quad I(v,f;P) \ge 2\, \rho(P,5)\ge 2\,\rho(P,4)\\
&I(g,C;P)\ge \min(I(v,C;P),I(u^r,C;P))\ge 2\,\rho(P,k)
\end{eqnarray*}
The equality $ I(v,f;P)=2\rho(P,5)$ holds only if $(C,P)\cong E_6$ and
$g_1(P)\ne 0$.

Assume $k=4$.  Then 
by the above inequality, we have 
$I(g,f;P)\ge 2\, \rho(4)$.
 The equality $I(g,f;P)=2\rho(P,4)$ holds if and only if  $(C,P)\cong A_j$
 singularity and $r=\rho(P,4)$.

Now we assume that $k=5$. Let $g$ be a conic
in $\Ker\, \si_{5,P}$,  written as 
$g=v\, g_1+u^r\, g_2$, $g_1,g_2\in \cO_P$.
By the above inequalities, we have
$I(v,f;P)\ge 2\, \rho(P,5)$ and the equality holds only for $E_6$.
The equality $I(g,f;P)=I(u^r,f;P)=2\rho(P,5)$ 
takes place if and only if  $(C,P)=A_j$ and
$r=\rho(P,5)$. Similarly the equality $I(g,f;P)=I(v,f;P)=2\rho(P,5)$
if  $(C,P)=E_6$ and 
$g_1(O)\ne 0$. In the last case,
 $g=0$ is smooth at $P$ and its
tangent line is the tangent cone of $(C,P)$.

The assertion II is obvious by the description of the ideal $\cJ_{P,k,6}$.
\qed
\begin{Remark}
Assume that $C$ is of torus type defined by $f_2^3+f_3^2=0$
 and $O\in \{f_2=f_3=0\}$.
Assume that $(C,O)\cong A_{3r-1}$. Then the cubic $f_3=0$ is smooth at
 $O$
and $I(f_3,f_2;O)=r$.
We may assume that $v=f_3$ and 
$f_2\equiv a\,u^r+\text{(higher terms)}$ modulo
$(v)$ for some $a\in \bfC^*$. Thus  $f_2\in \Ker\,\si_{5,O}$.
Assume next that $(C,O)$
is an $E_6$-singularity.
Then the cubic $f_3=0$ has a node at $O$ and $f_2$ is smooth at $O$ and 
the tangent cone is equal to  the tangent space of $f_2=0$
and we may assume that $f_2=v$
(\cite{Pho}). 
Thus again  $f_2\in \Ker\,\si_{5,O}$. Thus $f_2\in \si_{5,O}$ if $O$ is
 an 
simple inner singularity.
\end{Remark}

\begin{Lemma}\label{colinear lemma}
Let $\ell$ be a projective line and let  $P_1,P_2, P_3$ be  three points
 on $\ell$
(including the infinitely near point cases $P_1=P_2\ne P_3$ or
 $P_1=P_2=P_3$). 
Let $g=0$ be a conic on
 $\bfP^2$ which passes through them. 
Then $\ell$ divides $g$.
\end{Lemma}
{\em Proof.}
The assumption for the degenerated cases
$P_1=P_2\ne P_3$ or
 $P_1=P_2=P_3$
 implies
$I(g,\ell;P_1)=2$ or $3$ respectively.
Thus in any case, $I(g,\ell)\ge 3$. Thus by Bezout theorem,
this implies that $\ell\,|\,g$.\qed
\begin{Definition}
Let $P_1,\dots, P_\nu,\, \nu \le 3$ be $\rho(5)$-essential simple singularities on a
 reduced sextic $C$
such that $\sum_{i=1}^\nu \rho(P_i,5)\ge 3$.
We say that the singular points 
$\{(C,P_i);i=1,\dots,\nu\}$ are colinear  if there is a
projective line $L$ such that $P_i\in L$ for $i=1,\dots, \nu$ and
one of the following conditions is satisfied.

(1) $\nu=3$ or
(2) $\nu=2$, $I(L,C;P_1)\ge 4$ and $\rho(P_1,5)\ge 2$ and  $L$ is the principal tangential
 line at $P_1$  or 
(3) $\nu=1$ and $L$  is the principal tangential
 line at $P_1$, $I(L,C;P)\ge \, 6 $ and $\rho(P_1,5)\ge 3$.
\end{Definition}
Note that in either case, we have $\sum_{i=1}^\nu I(L,C;P_i)=6$ by
 Bezout theorem.

We  recall the definition of {\em sextics of linear type}.
Assume that $C$ is a sextics with $3A_5$ or $A_{11}+A_5$ or $A_{17}$.
We say that $C$ is {\em of linear type} if
there is a line $L\subset \bfP^2$ such that 
\[
L\cap C=\begin{cases}
\{P_1,P_2,P_3\},\, I(C,L;P_i)=2,\quad  & (C,P_i)\cong A_5,\,i=1,2,3\\
\{P_1,P_2\},\,I(C,L; P_1)=2,\, I(C,L;P_2)=4,\quad &(C,P_1)\cong A_5,\,
(C,P_2)\cong A_{11}\\
\{P_1\},\,I(C,L;P_1)=6,\quad &(C,P_1)\cong A_{17}
\end{cases}
\]
Let $\ell$ be the linear form defining $L$. Then we have shown in 
\cite{Reduced} that $C$ is
of torus type $f_3^2+f_2^3=0$ with $f_2=\ell^2$. Such a torus curve is
called
{\em of linear torus type}.
For further detail, see \cite{Reduced}.
Now we state our main theorem for sextics with simple singularities.
\begin{Theorem}\label{Global-rho} Assume that $C$ is a reduced
 sextics with only
 simple singularities.

\nin
{\rm (A)} The  homomorphism
$\si_4:H^0(\bfP^2,\cO(1))\to \bigoplus_{P\in \Si(C)}\cO_p/\cJ_{P,4,6}$
 is described as follows.
\nl
{\rm (a)} 
$\si_4$ is injective if $\rho(4)\ge 4$.  The case $\rho(4)> 4$
 does not exist.
\nl
{\rm (b)} If $\rho(4)\le 3$, the homomorphism
$\si_4:H^0(\bfP^2,\cO(1))\to \bigoplus_{P\in \Si(C)}\cO_p/\cJ_{P,4,6}$
 is surjective except the case \nl
\indent
 $(\sharp4):\, \rho(4)=3$ and  $C$ is of linear torus type.
\nl
In this case, $\dim \Coker\,\si_4=1$.

\nin
{\rm (B)} The  homomorphism
$\si_5:H^0(\bfP^2,\cO(2))\to \bigoplus_{P\in \Si(C)}\cO_p/\cJ_{P,5,6}$
 is
described as follows.
\nl
{\rm (c)}  $\si_5$
is surjective if $\rho(5)\le 5$.
\nl
{\rm (d)} If $\rho(5)\ge 6$, $\si_5$ is injective except the case that 
\nl\indent
$(\sharp5):\, C$ is  of torus type and $\rho(5)=6$. 
\nl
In the last case, 
 $\dim \Coker\,\si_5=1$.
\end{Theorem}

The proof of Theorem \ref{Global-rho} occupies the rest of this section.
Let $P_1,\dots,P_\nu$ be the 
$\rho(k)$-essential singularities.

{\em Proof of the assertion (A) in Theorem\ref{Global-rho}.}
We first prove the assertion (A).
By the classification tables
\cite{Pho,Oka-Pho2,Reduced} and Bezout theorem,  we observe that  $\rho(4)\le 4$ and
there exists a unique configuration $\Si=[3A_5,D_4]$ which satisfies the 
equality. Note that the rank is maximal for this configuration.

If there is  a line $\ell\in \Ker\, \si_4$, we have
$6\ge I(\ell,C)\ge 2\rho(4)$. This implies $\rho(4)\le 3$.
Thus assume that $\rho(4)\le 3$. The condition for a linear form $\ell$ to be in the
kernel
of $\si_4$ is given by $\rho(4)$ linear equations. If $\rho(4)<3$, they
are 
 independent by a direct computation. Thus $\dim\, \Ker \,\si_4= (3-\rho(4))$
and the surjectivity follows.
Assume that $\rho(4)=3$ and $\ell\in \Ker\, \si_4$.
appenWe may assume that $\ell$ is not a component of $C$. For the proof of this assertion, see
Appendix 1.
By Corollary \,\ref{Local-rho}, we have $I(\ell,C)\ge 6$. Thus we must
have the equality $I(\ell,C;P)=2\rho(P)$.
Thus  the singularities $P\in C$ with $\rho(P,4)>0$ must be $A_j,\, j\ge
5$ and they are on $\ell$ so that
 $I(\ell,C;P)=2\rho(P,4)$. As $\ell$ is unique up to a
multiplication 
of a constant by the property $I(\ell,C)\ge 6$, the assertion (b)
reduces to 

\begin{Lemma}\label{linear-torus-lemma}  Assume that $C$ is a sextic and  $\ell$ is a line and  one of the conditions are 
satisfied.
\nl
(i) $C\cap \ell=\{P_1,P_2,P_3\}$, $(C,P_i)\cong A_{j(i)}$ with
 $\rho(P_i,4)=1$, $5\le j(i)<11$ and $I(\ell,C;P_i)=2$ for $i=1,2,3$.
\nl
(ii) $C\cap \ell=\{P_1,P_2\}$, $(C,P_i)\cong A_{j(i)}$ with
 $\rho(P_1,4)=2$, $\rho(P_2,4)=1$,
$11\le j(1)<17,\, 5\le j(2)<11$ and $I(\ell,C;P_1)=4$ and $I(\ell,C;P_2)=2$.
\nl
(iii) $C\cap \ell=\{P_1\}$, $(C,P_1)\cong A_{j(1)}$ with
 $\rho(P_1,4)=3,\, 17\le j(1)<23$ and $I(\ell,C;P_1)=6$.
Under the above assumption, we conclude that
 $j(i)=6$ or $12$ or $18$, namely
the singularities on $C\cap \ell$ are   $3\,A_5$ 
or $A_{11}+A_5$ or $A_{17}$ respectively.
\end{Lemma} 
{\em Proof of Lemma\, \ref{linear-torus-lemma}.} Assume  that $r_i:=\rho(P_i,4)$ and 
$(C,P_i)=A_{j(i)}$ with $6r_i-1\le {j(i)}< 6r_i+5$.
Take  normal coordinates $(u,v)$ so that 
$C$ is defined by $v^2+u^{{j(i)}+1}+\text{(higher terms)}=0$.
The assumption $I(\ell,C;P_i)=2\, r_i$ implies that 
$\ell=a \, v+b\, u^{r_i}+\text{(higher
terms)},\, a,\, b\in \bfC$ and $b\ne 0$.
Let us consider the family of sextics
$C_t=f_t\inv(0)$, $f_t:=f(x,y)+t\, \ell(x,y)^6$. It is easy to observe that
$(C_t,P_i)\cong A_{6r_i-1}$ for a generic $\tau\in \bfC$. Thus for such a generic $\tau$,
$C_\tau$ is of linear type. Thus by Proposition 7 of \cite{Reduced} or by \cite{Tokunaga},
$C_\tau$ is of linear torus type written as 
$f_3(x,y)^2+c\,\ell(x,y)^6=0$ for some $c\ne 0$.
  This implies again $C=C_0$ is also of
linear torus type, as  $f_0=f_3^2+(c-\tau)\ell^6$.
Now the simple singularities of sextics of linear torus type are
$3\,A_5$, $A_{11}+A_5$ or $A_{17}$ (\cite{Reduced}).\qed


{\em Proof of the assertion (B) in Theorem \ref{Global-rho}}.
Now we consider the assertion (B).
Recall that $\dim H^0(\bfP^2,\cO(2))=6$.
Assume that $\rho(5)<6$. We will show the surjectivity of $\si_5$.
Let $P_1,\dots,P_\nu$ be the 
$\rho$-essential singularities and let  $r_i:=\rho(P_i,5)$.
For a conic $g$, the condition for $g$ to be in the kernel
of $\si(k,P)$ is given by $\rho(P,k)$ linear conditions.
Namely, taking a normal coordinate system
$(u_i,v_i)$ of the singularity $(C,P_i)$, the condition $g\in
\Ker\,\si_{5,P_i}$ is
given by (considering $g$ as a function of $u_i,\, v_i$)
$ g(0,0)=\frac{\partial^j g}{(\partial u_i)^j}(0,0)=0,\, j=1,\dots,
r_i-1$
and $i=1,\dots,\nu$.
Thus we have to show that these  conditions are independent
so that 
$\dim\,\Image\,\si_5=6-\dim \,\Ker\, \si_5=\rho(5)$.
If $\rho(5)\le 3$, the condition can be checked easily by case-by-case computation.
We consider the case $\rho(5)=4$. 

\nin
{\bf Colinear Case with $\rho(5)=4$}. Assume that $\{(C,P_i);i=1,\dots,\nu\}$ contains
 a colinear subset, say
$\{(C,P_i); i=1,\dots,\mu\},\, \mu\le \nu$. Let $L$ be a projective line 
such that $\sum_{i=1}^\mu I(L,C;P_i)=6$.
Then for any conic  $g\in \Ker\, \si_5$, we assert $L\, |\, g$.
Otherwise we have a contradiction:
\nl
 $2\ge \sum_{i=1}^\mu I(L,g;P_i)\ge 3$.
Thus we can write $g=L\times \ell$ with  $\ell=ax+by+c$.
Now the condition $g$ to be in the kernel of $\si_5$
is given by one linear condition among $a,\,b,\, c$ and
therefore the dimension of the kernel $\si_5$ is 2.
For example, consider the easiest case: $\nu=4$ and $\mu=3$.
Then $g$ is in  $\Ker \,\si_5$ if and only if $P_4\in \ell$.
The other possibilities are the cases $\nu=3,\, 2, \,1$.
Consider the case 
 $\nu=\mu=2$ and $r_1=r_2=2$. We may assume that  $L$ is 
the principal tangent direction of $(C,P_1)$ and $P_2$ is on $ L$.
Note that  the principal tangential direction of $(C,P_2)$ is
different from $L$.
Then the condition to be asked is
$P_2\in \ell$ as $I(g,C;P_2)\ge 4$.
We omit the proof of  the other cases as it is similar.

\nin
{\bf General case with $\rho(5)=4$.}
 Now we consider the generic case (non-colinear singularities) with 
$\rho(5)=4$.
The case $\nu=3$ or $4$ or $\nu=2$ and $r_1=r_2=2$ is easy to be
checked, by a direct computation
 after putting $P_i$ to
 suitable positions, say for example\nl
$P_1=(0,0)$, $P_2=(1,0)$, $P_3=(0,1)$ and $P_4=(1,1)$ for $\nu=4$
or\nl
$P_1=(0,0)$ with the principal tangential line $x=0$, $P_2=(1,1)$, 
$P_3=(1,-1)$ for $\nu=3$ and $r_1=2$, etc, using  $PGL(3,\bfC)$-action.

We consider two cases which is slightly less obvious:   (1) $\nu=1,\,
r_1=4$
or (2) $\nu=2,\, r_1=3$.
Assume that $(u,v)$ is a fixed coordinate system
which gives the normal form of $(C,P_1)$. We can choose affine $(x,y)$
coordinates
so that 
$ x=x(u,v)$, $y=y(u,v)$ where 
\[
x(u,v)=u+\sum_{i+j\ge 2} a_{i,j}\,u^i\, v^j,\quad
 y(u,v)=v+\sum_{i+j\ge 2} b_{i,j}\,u^i\, v^j
\]
Put $P_2=(\al,\be)$ in the case (2) $\nu=2$.
The genericity condition  implies that the principal tangential
line $L_v: \,y=0$ satisfies $I(L_v,C;P_1)\le 5$, or equivalently
  $b_{2,0}\ne 0$  for the case 1 (respectively $P_2\notin L$, so
$\be\ne 0$ in case 2).
Let $g=\sum_{i+j\le 2} c_{i,j} x^i\, y^j$ be a conic.
We consider the condition $g\in \Ker\, \si_{5,P_1}$ with $P_1=(0,0)$.
First we need $c_{0,0}=0$.
Let $g_1(u)=g(x(u,0),y(u,0))$.
By an easy computation, we get 
\[
 \Coeff(g_1,u)=c_{1,0},\,\Coeff(g_1,u^2)=c_{0,1} b_{2,0}+c_{2,0},
\, \Coeff(g_1,u^3)=c_{0,1} b_{3,0}+c_{1,1} b_{2,0}+2 c_{2,0} a_{2,0}
\]
Thus in the case 1, we can solve the equation
$g_1\equiv 0$ modulo $( u^4)$   
as 
\[\mathrm{g}(x,y)=
\mathit{c_{0,2}}\,y^{2} + ({\displaystyle \frac {\mathit{c_{0,1}}\,( - 
\mathit{b_{3,0}} + 2\,\mathit{b_{2,0}}\,\mathit{a_{2,0}})\,x}{\mathit{b_{2,0}}}} 
 + \mathit{c_{0,1}})\,y - \mathit{c_{0,1}}\,\mathit{b_{2,0}}\,x^{2}
\]
In the case 2,
the condition $g\in \Ker\, \si_5$ is equivalent to 
 $g(\al,\be)=0$  and $g_1\equiv 0$ modulo $(u^3)$  which gives 
\[\mathrm{g}(x,y)=
 - {\displaystyle \frac {(\beta \,\mathit{c_{1,1}}\,\alpha  + \beta 
\,\mathit{c_{0,1}} - \mathit{c_{0,1}}\,\mathit{b_{2,0}}\,\alpha ^{2})\,y^{2}
}{\beta ^{2}}}  + (\mathit{c_{1,1}}\,x + \mathit{c_{0,1}})\,y - \mathit{
c_{0,1}}\,\mathit{b_{2,0}}\,x^{2}
\]
In any  case, we see that the kernel is two-dimensional.

\vspace{.2cm}
Now we consider the case $\rho(5)\ge 5$ and we first show that
\begin{Claim} Assume that  $\rho(5)\ge 5$,   The kernel of $\si_5$ is 
at most
one-dimensional.
\end{Claim}
{\em Proof.}
Assume first that $\Ker\, \si_5$ contains two mutually coprime conics
$g,g'$.  This gives an obvious contradiction
\[
 4\ge I(g,g')\ge \sum_{i=1}^\nu I(g,g';P_i)\ge \sum_{i=1}^\nu r_i\ge 5
\]
Assume next that the kernel contains two linearly independent conics
$g,g'$ which has a common linear factor $\ell$.
Thus we can write
$g=\ell\ell_1$ and $g'=\ell \ell_1'$.
We may assume that $\ell$ is not a component of $C$. See Appendix 2.
 Put $P:=\ell_1\cap \ell_1'$. Then for any line $\xi$ through $P$, 
$\ell \xi$ is again in $\Ker\,\si_5$. 
\nl  (1) Assume that   $P\ne P_i,\, i=1,\dots,\nu$.
Then taking a generic $\xi$, the 
assumption $\ell\xi\in \Ker\, \si_5$ implies that $P_i\in \ell$
for $i=1,\dots,\nu$.
Then  we have a contradiction
\[
10\le 2\sum_{i=1}^\nu r_i \le \sum_{i=1}^\nu I(\ell\xi,C;P_i)=
I(\ell,C;P_i)\le 6
 \]
(2) Assume that $P=P_1$ for example. The $P_i\in \ell$ for $i=2,\dots,\nu$.
 By the assumption, we may assume that $\xi$ is a generic line through
$P_1$.
Thus $I(\xi,C; P_1)\le 3$. (Recall that the multiplicity of a simple
singularity is less than or equal to 3.) 
By Corollary \ref{Local-rho}, we have $\sum_{i=1}^\nu I(g,C;P_i)\ge 10$.
This gives  a contradiction
that 
\[
 \sum_{i=1}^\nu I(\ell,C;P_i)=\sum_{i=1}^\nu I(\ell\xi,C;P_i)
-I(\xi,C; P_1)\ge 10-3=7
\]

Now the surjectivity  of $\si_5$ for $\rho(5)\le 5$ is completed.

Assume that $\rho(5)\ge 6$. Assume that there is a conic $g\in \Ker\,
\si_5$. 
We assume  that $g$ is not a  component
of $C$. See Appendix 3 below for the proof.
Then by Corollary \ref{Local-rho}, we get 
$12\ge \sum_{i=1}^\nu I(g,C;P_i)\ge 2\, \rho(5)$
which implies $\rho(5)= 6$.
Thus we have $\sum_{i=1}^\nu I(g,C;P_i)=12$
and the conic $g=0$ intersects $C$ only at  $\rho$-essential singularities.
By a result of Tokunaga \cite{Tokunaga-torus}, $C$ is a sextics of torus
type. 
Thus we have  proved that $\si_5$ is injective 
if $\rho(5)\ge 7$ or $\rho(5)=6$ and $C$ is not of torus type.
This proves the assertion (B) of Theorem \ref{Global-rho}.
\qed

If $C$ has the rank 19, by Proposition \ref{simple-rho} and 
 by the classification table in 
\cite{Yang}, we see that 
$C$ is of torus type and  $C$ has maximal rank 19. This  suggests us:

\nin
\begin{Conjecture}\label{conjecture}
 Assume that $C$ is a reduced sextic and
assume that $\rho(5)\ge 7$. Then $C$ is a sextic of torus type.
\end{Conjecture}

\begin{Corollary}\label{Simple-list} Assume that $C$ is a reduced sextic
 of torus curve with only simple singularities.
Then \begin{enumerate}
\item $\rho(5)=9$ if and only if $C$ has 9 cusps and in this case
  $ \De(t)=(t^2-t+1)^3$.
\item Assume that  $\rho(5)=8$.
\begin{enumerate}
\item $\rho(4)\le 2$ if and only if 
$\Si(C)$ is one of the following.
\[
[8A_2],\,  [8A_2,A_1],\,[6A_2,E_6],\,[6A_2,A_5],
\,[6A_2,A_5,A_1], \,[4A_2,2A_5],\,[4A_2,A_5,E_6]
\]
In this case, $C$ is irreducible and   $ \De(t)=(t^2-t+1)^2$.
\item $\rho(4)\ge 3$ if and only if
$\Si(C)=[3A_5,2A_2]$. In this case, $C$ is of linear torus type and
$C$      consists of two cuspidal cubic components and 
$\tilde\De(t)=(t^2-t+1)^2 (t^2+t+1)$.\end{enumerate} 
\item The case $\rho(5)=7$ and $\rho(4)=4$ 
occurs if and only if $\Si(C)=3A_5+D_4$.
In this case, $\tilde \De(t)=(t^2-t+1) (t^2+t+1)$.
There are three components in the moduli space:
\nl
{\rm (a)} $C=B_4+ B_1+B_1'$, where $B_4$ is a quartic and
$B_1,\, B_1'$ are  two flex tangents of $B_4$ meeting at a point on
      $B_4$
 and $\Si(C)=[3A_5,D_4]$ .\nl
{\rm (b)} $C$ has  three conic components so that $\Si(C)=[3A_5,D_4]$.
\nl
{\rm (c)} $C$ is of linear torus with $\Si(C)=[3A_5,D_4]$. In this case,
$C$ consists of a cubic component and three flex tangent lines.\nl
\item Assume that 
$\rho(5)\le 7$ and $\rho(4)=3$ and $C$ is of linear torus type.
Then  $\tilde\De(t)=(t^2-t+1)(t^2+t+1)$.
\item For other sextics of  torus type 
which  is not of linear torus type, we have  $\tilde \De(t)=t^2-t+1$.
\end{enumerate}
The irreducible sextics of torus type can occur only in 1, 2-(a)  or  5.
\end{Corollary}
The proof of Corollary is immediate from Theorem \ref{Global-rho} and 
 the classification tables in \cite{Pho,Oka-Pho2,Reduced}.
\subsection{Sextics of non-torus type}
We consider irreducible sextics of non-torus type with only simple 
singularities. 
We consider the factor $t^2+t+1$. This is determined by $\rho(4)$.
If $\rho(4)\le 3$ and $C$ is not of linear torus type, we know by
Theorem \ref{Global-rho} that the factor $t^2+t+1$ does not appear.
Now we consider the possibility that $\rho(4)\ge 4$.
Recall that   the total Milnor number  is bounded by 19.
Put $\hat\mu(4)=\sum_{\rho(P,4)>0}\mu(C,P)$.
Proposition \ref{simple-rho} implies that $\rho(4)\le 4$ and the 
possibility of $\rho(4)=4$ is the case that $\Si(C)$ contains
one of the following.
\begin{enumerate}
\item $\hat\mu(4)=16$: $4D_4\subset \Si(C)$.
\item $\hat\mu(4)=17$: $3D_4+X_5\subset \Si(C),\, X=D,\,  A$.
\item $\hat\mu(4)=18$: $3D_4+X_6\subset \Si(C),\,X=D,\, A,\, E$ or
$2D_4+X_5+Y_5\subset \Si(C),\, X,Y=D,\, A$ or $D_{10}+2D_4\subset \Si(C)$.
\item $\hat\mu(4)=19$: $3D_4+X_7\subset \Si(C),\,X=D,\, A,\, E$
or $2D_4+X_j+Y_k\subset \Si(C),\, k+j=11,\, k,j\ge 5, \, X=D,\, A,\, E$
or $2D_4+X_{11}\subset \Si(C),\, X=D$,  $A$,
or $D_4+X_5+Y_5+Z_5\subset \Si(C)$, $X,Y,Z=$  $D,A$
or $D_4+D_{10}+X_5\subset \Si(C)$ with $X=D$, $A$.
\end{enumerate}
On the other hand, the irreducibility gives the Pl\"ucker inequality:
\[
\sum_{\rho(P,4)>0} \de(C,P)\le 10
\]
As  $\de(C,P)=(2\mu(C,P)+r-1)/2$ where $r$ is the number of irreducible
components
at $P$ ( so $\de(D_4)=\de(A_5)=3$ and so on), 
we can easily see that none of the
above
configurations  satisfy the Pl\"ucker inequality.
Thus  the factor $t^2+t+1$ does not appear in 
Alexander polynomials of irreducible sextics of non-torus type with
simple singularities. By Theorem \ref{Global-rho} and Corollary \ref{Simple-list} we obtain:
\begin{Corollary}\label{irreducible-simple}
Assume that $C$ is an irreducible sextics with only simple singularities.
Then the Alexander polynomial $\De(t)$ takes the form
$(t^2-t+1)^\al,\, \al=0,1,2,3.$
\end{Corollary}
\subsection{Appendix}
{\bf  1.} We prove that {\em there does not exist a sextics with
simple singularities
and $\rho(4)\ge 3$ which has  a line component $\ell=0$ such that $\ell\in
\Ker\,\si_4$.}

Assume that there exists a sextic $C=L\cup B_5$ with
a line component $L=\{\ell=0\}$ such that $\ell\in \Ker\,\si_4$.
Put $P_1,\dots,P_\nu$ be $\rho(4)$-essential singularities. We assume 
that $(C,P_i)$ are simple singularities. We will show that $\rho(4)\le 2$.
By the
assumption $\ell\in \Ker\,\si_4$ implies that $P_i\in L$ for
$i=1,\dots, \nu$.  Note also
$\sum_{i=1}^\nu I(L,B_5;P_i)\le 5$. Thus $\nu\le 2$.
\nl
(1) If  $(C,P_i)\cong A_{2k-1}: f(u,v)= v^2-u^{2k}=0$, we may assume 
 $\ell=v-u^k$, $I(L,B_5;P_i)=k\le 5$ and 
$ \rho(4,P)= 1$  for $ k=3,4,5$.
Actually, $f(u,v)$ may have higher terms 
and
then $\ell=v-u^k+\text{(higher terms)}$.
We ignored these higher terms also in th following cases
 but the computation are the same.
\nl
(2) If  $(C,P_i)\cong E_7: v^3+vu^3=0$, 
$\ell=v$, 
 $I(L,B_5;P_i)=3$ and $\ell\in \Ker\,\si_{5,P_i}$ and $\rho(4,P_i)=1$.
\nl
(3) If $(C,P_i)\cong D_k: v^2u-u^{k-1}=0$, and  $\ell=\{u=0\}$. Note
that  $\ell\in
\Ker\,\si_{4,P_i}$ if and only if
$k\le 9$. In this case,
 $I(L,B_5;P_i)=2$ and and $\rho(4,P_i)=1$.
For $k=2m$, $D_{2m}$ has three components and if
 $\ell=\{v-u^{m-1}=0\}$,
$I(L,B_5;P_i)=m\le 5$.  Thus $\rho(4,P_i)=1$ for $m=2,3,4$
and $\rho(4,P_i)=2$ for $m=5$.

Thus in any case, we have observed that
$ \rho(4,P_i)\le \frac 12 I(\ell,B_5;P_i)$.
This implies that
\[
 \rho(4)=\sum_{i=1}^\nu \rho(4,P_i)\le \sum_{i=1}^\nu I(\ell,B_5;P_i)\le
 \frac 52
\]

{\bf  2.}
Now we show that {\em if $\ell$ is a component of $C$
and if $\ell\ell_1,\,\ell\ell_1'$ are 
linearly independent conics in $ \Ker\,\si_5$, then $\rho(5)\le 4$.}

For sextics of torus type, the assertion follows from the list of 
configurations with a line component in \cite{Reduced}.
Assume that $C$ need not be of torus type and  $\ell$ is a component
of $C$ and $P_1,\dots, P_\nu$ are on $\ell$.
Let $C=\ell+B_5$ where  $B_5$ is a reduced curve of degree 5.
The possible reducible singularity
 $(C,P_i)$ is either $A_{2k-1}$ or $E_7$ or  $D_{j}$
and $\ell$ must be a smooth component. Let $P=\ell\cap \ell_1$.
We can see as before that $P_i\in \ell$ and $\nu\le 2$
as $\sum_{i=1}^\nu I(\ell,B_5;P_i)\le 5$.
\nl
(1) For $A_{2k-1}: v^2-u^{2k}=0$, if $\ell=v-u^k$, $I(v-u^k,B_5;O)=k\le 5$ and 
$\ell\in \Ker\,\si_{5,P_i}$ and $\rho(5,P_i)\le [2k/3]$.
\nl
(2) For $E_7: v^3+v\,u^3=0,\, s\ge 5$, $\ell=v$, 
 $I(L,B_5;O)=3$ and $\ell\in
 \Ker\,\si_{5,P_i}$ and $\rho(5,P_i)=2$.
\nl
(3) For $D_k: v^2u-u^{k-1}=0$, 

if (3-1) $\ell=\{u=0\}$,
 $I(L,B_5;O)=2$ and $\rho(5,P_i)=[k/3]$. The assumption $\ell\in
 \Ker\,\si_{5,P_i}$ implies
 $k=4,\,5$.

 For $k=2m$ and
if (3-2) $\ell=\{v-u^{m-1}=0\}$,
$I(L,B_5;O)=m$ and $\ell\in \Ker\, \si_{5,P_i}$.
As $I(L,B_5;P_i)\le 5$, we need $m\le 5$.

We observe that in any case, 
\begin{eqnarray}\label{rho-estimate}
 \ell\in \Ker\,\si_{5,P_i}\implies \rho(C,P_i)\le [\frac 23 I(\ell,B_5;P_i)]
\end{eqnarray}
In the above discussion, the higher terms
in the defining equations are ignored but the same assertion holds.
Put $P=\ell_1\cap \ell_1'$. We can assume that $\ell_1$ is a generic
line through $P$. First, assume that $P$ is not an
$\rho(5)$-essential singularity.
This implies that $\ell\in \Ker\,\si_5$. Then  we have
$\rho(5)\le \sum_{i=1}^\nu \frac 23 I(\ell,B_5;P_i)\le \frac{10}3$ by (\ref{rho-estimate}).
Next we may assume that $P=P_1$, $(C,P_1)\cong D_k$ and $\ell=u$.
The assumption 
$\ell\ell_1\in \Ker\,\si_{5,P_1}$ is equivalent to $u^2\in \Ker\,
\si_{5,P_1}$.
Thus $\rho(C,P_1)=2$.
The possibilities are
 $\nu=2$: $D_i+A_{2j-1},\, (i=6,7,8,\, j=2,3)$ or $D_i+E_7,\, (i=6,7,8)$
or 
$\nu=1$ and  $D_i,(i=6,7,8)$.
Thus in any case we have $\rho(5)\le 3$.

{\bf  3.} Assume that $C$ is a reduced sextic with only
simple singularities.
We show that {\em if $g=0$ a component of a sextics $C$ 
and  $g\in \Ker\,\si_5$, then
 $\rho(5)\le 5$.}

Assume that $C=C^2\cup C^4$ where  $C^2$ is defined by  $g=0$
with $g\in \Ker\, \si_5$.
Let $P_1,\dots,P_\nu$ be $\rho$-essential singularities.
Then  $P_i\in C^2$ for  $i=1,\dots, \nu$.
The singularity $(C,P_i)$ is reducible and thus it is one of 
$A_{2j-1}$, $E_7$ or $D_j$. Note that $(\ref{rho-estimate})$ is  true
under the substitution $B_5\to C^4$.
Thus we get 
\[
\rho(5)= \sum_{i=1}^\nu \rho(5,P_i)\le \frac 23 \sum_{i=1}^\nu
I(C^2,C^4;P_i)\le 
\frac{16}3
\]
which implies $\rho(5)\le 5$.
\section{Non-simple singularities}
In this section, we consider sextics of torus type with some non-simple
 singularities.
By 
Pho \cite{Pho}, non-simple singularities on sextics of
torus type are
$B_{3,2k}\,(k=2,\dots,6)$, $B_{4,6},B_{6,6}$, $C_{3,k}\,(k=7,8,9,12,15)$,
$C_{6,3k}\,(k=2,3,4)$, $Sp_1, \,Sp_2$ and $D_{4,7}$ where
\begin{eqnarray*} 
B_{m,n}:&\quad f(u,v)=v^m+u^n+\text{(higher terms)}=0\qquad {\rm \text{(Pham-Brieskorn type)} }\label{bmn-def}\\ 
C_{m,n}:&\quad f(u,v)=v^m+u^2v^2+u^n+\text{(higher terms)}=0,\frac 2m+\frac 2n
< 1\label{cmn-def}\\ 
D_{4,7}:&\quad v^4+u^3 v^2\,+\, a u^5 v\,+\, bu^7=0,\qquad a^2-4b\ne 0\\
Sp_1:&\quad f(u,v)=f_2(u,v)^3+f_3(u,v)^2,~ 
\begin{cases} 
f_3(u,v)=&v^2- u^3+\text{(higer terms)}\\ 
f_2(u,v)=&cuv+\text{(higer terms)},~c\ne 0 
\end{cases} \\
Sp_2:& f_3(u,v)^2-\, c\,v^6=0,\, c\ne 0,\quad
f_3(u,v)=v^2- u^3+\text{(higer terms)}
\end{eqnarray*} 
Note that $B_{m,n}$ and $C_{m,n}$ have non-degenerate 
Newton boundaries and   
$Sp_1, Sp_2$  are  degenerated in the sense of Newton boundary (\cite{Okabook}).
The following describes the $\rho(k)$-invariants and the $\cJ_{P,k,6}$ 
ideals of the above non-simple singularities. 
 We use the notation $\rho_{3,5}(P)=(\rho(5,P),\rho(4,P),\rho(3,P))$.
\begin{Lemma}\label{non-simple-rho} Assume that $C$ is a sextics and $P$ is a singular point 
defined by the above equation. Then the ideal
$\cJ_{P,k,6}$ and $\rho(k,P)$ are given
 as follows.
\begin{enumerate}
\item $B_{3,6} $, $C_{3,7}$ and $C_{3,8}$:

 $\cJ_{P,k,6}=\langle u^3,uv, v^2\rangle,
\langle u^2,v\rangle,\,
\langle u,v\rangle$ and $\rho_{3,5}(P)=(4, 2, 1)$. 
\item
$C_{3,9}$ and $B_{3,8}$:
$\cJ_{P,k,6}=\langle u^4,v^2,u v\rangle,
\langle u^2,v\rangle,\,
\langle u,v\rangle$ and $\rho_{3,5}(P)=(5, 2, 1)$. 
\item $C_{6,6}$: $\cJ_{P,k,6}=\langle u^3,v^3,u v\rangle,
\langle u^2,v^2,uv\rangle,\,
\langle u,v\rangle$ and $\rho_{3,5}(P)=(5, 3, 1)$. 
\item $B_{4,6},\,D_{4,7}$: $\cJ_{P,k,6}=\langle u^3,v^2,vu^2\rangle,
\langle u^2,v^2,uv\rangle,\,
\langle u,v\rangle$ and $\rho_{3,5}(P)=(5, 3, 1)$. %
\item $C_{3,12}$: $\cJ_{P,k,6}=\langle u^5,v^2,u v\rangle,
\langle u^2,v\rangle,\,
\langle u,v\rangle$ and $\rho_{3,5}(P)=(6, 2, 1)$.
\item $B_{3,10}$:
$\cJ_{P,k,6}=\langle u^5,v^2,u v\rangle,
\langle u^3,v\rangle,\,
\langle u,v\rangle$ and $\rho_{3,5}(P)=(6, 3, 1)$.
\item $C_{6,9}$: $\cJ_{P,k,6}=\langle u^4,v^3,u v\rangle,
\langle u^2,v^2,uv\rangle,\,
\langle u,v\rangle$ and $\rho_{3,5}(P)=(6, 3, 1)$.
\item $Sp_1$:
$\cJ_{P,k,6}=\langle u^4,v^2-u^3,v^3,u^2 v\rangle,
\langle u^2,uv,v^2\rangle,\,
\langle u,v\rangle$ and $\rho_{3,5}(P)=(6, 3, 1)$.
\item $C_{3,15}$:
$\cJ_{P,k,6}=\langle u^6,v^2,u v\rangle,
\langle u^3,v\rangle,\,
\langle u,v\rangle$ and $\rho_{3,5}(P)=(7,3,1)$. 
\item $C_{9,9}$:
$\cJ_{P,k,6}=\langle u^4,v^4,u v\rangle,
\langle u^2,v^2,uv\rangle,\,
\langle u,v\rangle$ and $\rho_{3,5}(P)=(7, 3, 1)$. 
\item $C_{6,12}$:
$\cJ_{P,k,6}=\langle u^5,v^3,u v\rangle,
\langle u^3,v^2,uv\rangle,\,
\langle u,v\rangle$ and $\rho_{3,5}(P)=(7, 4, 1)$. 
\item $Sp_2$:  $\cJ_{P,k,6}=\langle u^4,v^2-u^3,v^3,u v^2\rangle,
\langle u^3,uv,v^2\rangle,\,
\langle u,v\rangle$ and $\rho_{3,5}(P)=(7, 4, 1)$.
\item $B_{3,12}$: $\cJ_{P,k,6}=\langle u^6,v^2\rangle,
\langle u^4,v\rangle,\,
\langle u^2,v\rangle$ and $\rho_{3,5}(P)=(8, 4, 2)$.
\item $B_{6,6}$:
$\cJ_{P,k,6}=\langle u^4,u^3v,u^2v^2,uv^3,v^4\rangle,
\langle u^3,u^2v,uv^2,v^3\rangle,\,
\langle u^2,uv,v^2\rangle,\,$
$\langle u,v\rangle$ and $\rho_{3,5}(P)=(10, 6, 3)$ and $\rho(2,P)=1$.
\end{enumerate}
In the above list, $\rho(P,2)=0$  except for $B_{6,6}$ and we omitted it.
\end{Lemma} 
{\em Proof.}
Observe that $\cJ_{P,5,6}$ for $Sp_1$ and $Sp_2$  are not  monomial ideals.  
Except for these two singularities, the assertion follows from Lemma
\ref{non-degenerate-rho}. 
We will show the computation for $C_{3,9}$ and leave the other cases to 
the reader.
Assume that $(C,P)\cong C_{3,9}$. Thus the defining equation is
\[
 C_{3,9}:\quad f(u,v)=v^3+u^2v^2+u^9=0
\]
We use Lemma \ref{non-degenerate-rho}. We have two faces corresponding
to  the face function $v^3+u^2v^2$ and $u^2v^2+u^9$ whose weight vectors
are
\[
 P_1=\begin{pmatrix} 1\\2\end{pmatrix},
\quad P_2=\begin{pmatrix} 2\\7\end{pmatrix},
\]
Let $k_i,m_i$ be the multiplicity of $\pi^*du\land dv$ and $\pi^* f$
along the divisor ${\hat E}(P_i)$.
Thus the ideal $\cJ_{P,k,6}$ is generated by 
the functions whose pull back have zeros of multiplicity at least
$[5m_1/6-k_1]$ and $[5m_2/6-k_2]$  along ${\hat E}(P_1)$ and
  ${\hat E}(P_2)$ respectively. These integers are given by 3 and 7 for $k=5$
and  
  $2$ and $4$ for  $k=4$ and $1$ and $1$ for $k=3$.
On the other hand,  $(\pi^*u)={\hat E}(P_1)+{\hat E}(P_2)+(others)$ and $(\pi^*v)=2{\hat E}(P_1)+7{\hat E}(P_2)+(others)$ .
Thus we can see easily that 
$\cJ_{P,k,6}$ is generated by the monomials
\[
 \cJ_{P,5,6}=\langle u^4,v^2,uv\rangle,
\,
\cJ_{P,4,6}=\langle u^2,v\rangle, \,
\cJ_{3,P,6}=\langle u,v\rangle
\]

For $Sp_1$ and $Sp_2$, we proceed twice toric modifications to obtain their
resolutions.

{\bf I.}
We show the assertion for the case $Sp_1: f(u,v)=(v^2-u^3)+c\,(uv)^3=0$.
We take the first toric modification $\pi_1: X_1\to \bfC^2$
 with respect to the regular simplicial
cone
with vertices
\[
 T_0=\begin{pmatrix} 1\\0\end{pmatrix},
T_1=\begin{pmatrix} 1\\1\end{pmatrix},
T_2=\begin{pmatrix} 2\\3\end{pmatrix},
T_3=\begin{pmatrix} 1\\2\end{pmatrix},
T_4=\begin{pmatrix} 0\\1\end{pmatrix}
\]
The weight vector $T_2$ corresponds to the unique edge of 
the Newton boundary
$\Ga(f;(u,v))$.
Consider the toric chart $\Cone(T_2,T_3)$ and denote its
toric coordinates by  $(u_1,v_1)$.
Then we have
\[
 u=u_1^2 v_1,\, v=u_1^3 v_1,\, \pi_1^*f(u_1,v_1)=u_1^{12}
 v_1^6((v_1-1)^2+c\,{v_1}^3 u_1^3)
\]
and the two form $K=du\land dv$ is  shifted as
\[
 du\land dv=u_1^4v_1^2 du_1\land dv_1
\]
The strict transform $C'$ of $C$ intersects with ${\hat E}(T_2)$.
Take a new coordinate system $(u_1,v_1')=(u_1,v_1-1)$.
Now the local behavior of $\pi_1^*f$ 
on ${\hat E}(T_2)$ is written as 
$\pi_1^*(u_1,v_1')=u_1^{12}({v_1'}^2+c\, u_1^3)+\text{(higher terms)}$.
Thus we proceed the second toric modification
$\pi_2:X_2\to X_1$ with respect to the exactly
same simplicial cone:
\[
S_0=\begin{pmatrix} 1\\0\end{pmatrix},
S_1=\begin{pmatrix} 1\\1\end{pmatrix},
S_2=\begin{pmatrix} 2\\3\end{pmatrix},
S_3=\begin{pmatrix} 1\\2\end{pmatrix},
S_4=\begin{pmatrix} 0\\1\end{pmatrix}
\]
Let $(u_2,v_2)$ be the toric coordinates with respect to 
$\Cone(S_2,S_3)$. Then 
$u_1=u_2^2 v_2,\, v_1'=u_2^3v_2^2$ and 
\[
u_1^{4} du_1\land dv_1'= u_2^{12}v_2^6 du_2\land dv_2
\]
Thus putting $\pi=\pi_1\circ \pi_2:X_2\to \bfC^2$, we have 
\begin{eqnarray*}
&(K)={\hat E}(T_1)+4{\hat E}(T_2)+2{\hat E}(T_3)      +5{\hat E}(S_1)+12 {\hat E}(S_2)+6{\hat E}(S_3)\\
&(\pi^*f)=4{\hat E}(T_1)+12{\hat E}(T_2)+6{\hat E}(T_3)+14{\hat E}(S_1)+30 {\hat E}(S_2)+ 15{\hat E}(S_3)\\ 
&(\pi^* u)={\hat E}(T_1)+2{\hat E}(T_2)+{\hat E}(T_3)+2{\hat E}(S_1)+4{\hat E}(S_2)+2{\hat E}(S_3)\\
&(\pi^* v)={\hat E}(T_1)+3{\hat E}(T_2)+ 2{\hat E}(T_3)+3{\hat E}(S_1)+6{\hat E}(S_2)+3{\hat E}(S_3)
\end{eqnarray*}
The divisors 
corresponding to $\cJ_{P,5,6}$, $\cJ_{P,4,6}$ and $\cJ_{P,3,6}$ are
given by
\begin{eqnarray*}
&\cJ_{P,5,6}=\{\phi; (\phi)\ge 2{\hat E}(T_1)+6{\hat E}(T_2)+3{\hat E}(T_3)+6{\hat E}(S_1)+13{\hat E}(S_2)+6{\hat E}(S_3)\}\\
&\cJ_{P,4,6}=\{\phi; (\phi)\ge {\hat E}(T_1)+4{\hat E}(T_2)+2{\hat E}(T_3)+4{\hat E}(S_1)+8{\hat E}(S_2)+4{\hat E}(S_3)\}\\
&\cJ_{P,3,6}=\{\phi;(\phi)\ge 
2{\hat E}(T_2)+2 {\hat E}(S_1)+3{\hat E}(S_2)+{\hat E}(S_3)\}
\end{eqnarray*}
Let us consider $\cJ_{P,5,6}$.
It is easy to see that $u^4,u^2v,uv^2, v^3\in \cJ_{P,5,6}$. Furthermore
we observe that $v^2-u^3$ is also in the ideal and they generate the
ideal. This implies that $\rho(P,5)=6$.
For $k=4,3$, it is easy to see that 
$ \cJ_{P,4,6}=\langle u^2,uv, v^2\rangle$ and $\cJ_{P,3,6}=\langle
u,v\rangle$. So  $\rho(P,4)=3$ and $\rho(P,3)=1$.
Observe that $\cJ_{P,5,6}$ is not a monomial ideal.

{\bf II.} Now we consider $Sp_2$ defined by 
$f(u,v)=(v^2-u^3)^2+ c\, v^6=0$. First we take the same toric modification
$\pi_1:X_1\to \bfC^2$
 with respect to the regular simplicial
cone
with vertices
\[
 T_0=\begin{pmatrix} 1\\0\end{pmatrix},
T_1=\begin{pmatrix} 1\\1\end{pmatrix},
T_2=\begin{pmatrix} 2\\3\end{pmatrix},
T_3=\begin{pmatrix} 1\\2\end{pmatrix},
T_4=\begin{pmatrix} 0\\1\end{pmatrix}
\]
Then with respect to the toric coordinate $(u_1,v_1)$
of  $\Cone(T_2,T_3)$,
the pull-back is written as
\begin{eqnarray*}
\pi^*f(u_1,v_1)&=u_1^{12}v_1^{6}\{(v_1-1)^2+c \,u_1^6v_1^2 \}\\
&=u_1^{12}\{{v_1'}^2+c\, u_1^6+\text{(higher terms)}\}
\end{eqnarray*}
Thus  we need one more toric modification $\pi_2:X_2\to X_1$
 with respect to the covectors
\[
R_0 = \begin{pmatrix} 1\\0\end{pmatrix},\,
R_1= \begin{pmatrix} 1\\1\end{pmatrix},\,
R_2= \begin{pmatrix} 1\\2\end{pmatrix},\,
R_3= \begin{pmatrix} 1\\3\end{pmatrix},\,
R_4=\begin{pmatrix} 0\\1\end{pmatrix}
\]
where the divisor $ {\hat E}(R_0)={\hat E}(T_2)$ and $ {\hat E}(R_3)$
corresponds to the
face ${v_1'}^2+c u_1^6$.
By a similar computation, we can show that 
$\cJ_{P,5,6}=\langle u^4, v^3,uv^2,u^3-v^2\rangle$
and $\cJ_{P,4,6}=\langle u^3,v^2,uv\rangle$ and 
$\cJ_{P,3,6}=\langle u,v\rangle$. This implies
$\rho(P,5)=7,\rho(P,4)=4$ and $\rho(3)=1$.
This completes the proof of Lemma \ref{non-simple-rho}.
\section{Alexander polynomial for sextics of torus type}
\begin{Lemma}\label{Local-rho-NS} 
Assume that $C$ is a sextic of torus type
and  $(C,P)$ is a non-simple singularity. 
Then  
$I(g,C;P)> 2 \rho(k,P)$  for 
 $g\in \cO_P\cap \cJ_{P,k,6}$.
\end{Lemma}
{\em Proof.}
The proof follows from
 the ideal $\cJ_{P,k,6}$ description given by Lemma \ref{non-simple-rho}. We observe that $\rho(5)\ge 7$ for sextics of torus type with
at least one non-simple singularity (see \cite{Pho,Oka-Pho2,Reduced}).
For example, consider the case $B_{3,6}\in C$ and assume that
$(C,P)\cong B_{3,6}$.  Then $\rho(P,5)=4$ and $\rho(P,4)=2$.
Assume that $g\in \cO_P\cap \cJ_{P,k,6}$.
For  $k=5$, we see
that 
$g(u,v)=u^3\, a(u,v)+v^2\,b(u,v)+\, uv\,c(u,v)$ with some $a,b,c\in
\cO_P$.
As $I(u^3,C;P)=9$, $I(v^2,C;P)=12$ and $I(uv,C;P)=9$, we get 
$I(g,C;P)\ge 9$. For $k=4$, we can write
$g(u,v)=u^2 a(u,v)+v\, b(u,v)$.
This implies that $I(g,C;P)\ge 6$.
 This proves the assertion. For other singularities,
similar argument using the linear combination of generators of
$\cJ_{P,5,6}$ and the normal form of the singularities proves the
assertion.
\qed

Now the assertion corresponding to Theorem \ref{Global-rho} 
takes the following form.
\begin{Theorem}\label{Global-rho-bis} Assume that $C$ is a reduced
 sextics of torus type with at least one non-simple singularities. 

\nin
{\rm (A)}
 The homomorphism
$\si_5:H^0(\bfP^2,\cO(2))\to \bigoplus_{P\in \Si(C)}\cO_p/\cJ_{P,5,6}$
is injective.

\nin
{\rm (B)}
{\rm (a)} The homomorphism
$\si_4:H^0(\bfP^2,\cO(1))\to \bigoplus_{P\in \Si(C)}\cO_p/\cJ_{P,4,6}$
 is injective for $\rho(4)\ge 4$.

{\rm (b)} If $\rho(4)\le 3$, the homomorphism
$\si_4:H^0(\bfP^2,\cO(1))\to \bigoplus_{P\in \Si(C)}\cO_p/\cJ_{P,4,6}$
 is surjective.
\end{Theorem}
{\em Proof.}
Let $P_1,\dots,P_\nu$ be the 
$\rho$-essential singularities.
We first prove the assertion (A).
Assume that $g$ is a conic in $\Ker\,\si_5$. By the classification tables in 
\cite{Pho,Oka-Pho2,Reduced}, we see that $\rho(5)\ge 7.$ Thus by
Corollary \ref{Local-rho} and Lemma \ref{Local-rho-NS}, we get
$I(g,C)\ge 2 \rho(5)\ge 14$, a contradiction to Bezout theorem if $g=0$
is not a components of $C$. 
Now assume that $g=0$ is a conic component $C^2$
of $C$. 
Then $C=C^2+C^4$ where $C^4$ is a reduced quartic.
Assume that $P_1$ is a non-simple singularity. Then  by th assumption
$g\in \Ker\, \si_5$,
 $P_2\in C^2\cap C^4$. 
By Lemma \ref{non-simple-rho}, we can see that no smooth component is in
$\cJ_{P_1,5,6}$. Thus $(C^2,P_1)\cong A_1$. In particular, $C^2$
is a union of two lines. 
The non-simple singularities which can  have two smooth components
 are
$C_{6,6}$, $C_{6,12}$ or  $B_{6,6}$.
Assume that 
$(C,P_1)\cong C_{6,6}$.
By the classification of reduced sextics with 
$C_{6,6}$ (\cite{Reduced}), possible configurations are
(a) $\Si(C)=C_{6,6}+2A_2+2A_1$  and $C$ has  two line components and a quartic $B_4$
or (b) $\Si(C)=C_{6,6}+A_5+2A_1$ and $C$ has two line components and two conic 
components. In the case of (a), two $A_2$ are on the quartic.
In case (b), $A_5$ has to be on the intersection of two conics. In any case, 
$g$ can not be in $\Ker\,\si_5$.
Assume that $(C,P_1)\cong C_{6,12}$. By \cite{Reduced}, there are no
possibility of sextics with two linear components.
Assume that $(C,P_1)\cong B_{6,6}$. This implies that $C$ consists of 6
lines meeting at $P_1$. By Lemma \ref{non-simple-rho}, no conic can be
contained in $\cJ_{P_1,5,6}$.
 Thus $g\notin \Ker\,\si_5$.

Now we prove the assertion (B). First we observe that 
$\rho(P,4)\ge 2$ for any non-simple singularities
which appears on sextics of torus type.
Assume that $\rho(4)=2$ and assume that $\ell,\ell'$ be independent
linear forms in $\Ker\,\si_4$.
Then this gives an contradiction $I(\ell,\ell')\ge 2$.
Now we assume that $\rho(4)\ge 3$ and we show that $\si_4$ is
injective.
First we assume that $\ell$  is a line which is not a component of $C$.
Then $\ell$ can not be in $\Ker\, \si_4$ as
otherwise we have a contradiction:
\[
6\ge  \sum_{P\in \ell}I(\ell,C;P) > 2 \rho(4)\ge 6
\]
Now we prove that $C$ does not have a line component $\ell$ which is
 in $\Ker\,\si_4$. Assume that $P$ is a non-simple singularity
and assume that $\ell=0$ is a line component such that 
$\ell\in \Ker\, \si_{P,4}$. Put $C=B_1+C^5$ where $B_1=\{\ell=0\}$ and 
$C^5$ is the union of other components.
The possibility for 
$(C,P)$ with $\ell\in \Ker\, \si_{P,4}$ is, by Lemma \ref{non-simple-rho},
 one of 
$B_{3,6},C_{3,k}, B_{3,10}, B_{3,12}$.
By \cite{Reduced}, $B_{3,10}$ does not appear on reduced sextics and 
$B_{3,12}$ is only possible for sextics with three conic components.
For $B_{3,6}$ or $C_{3,k}$, $\rho(P,4)=2$ and we must have other
 singularity $P'$ with $\rho(P',4)>0$.
However by \cite{Reduced}, we know that 
$I(B_1,C^5;P)=4$ or $5$.  Thus $B_1$ can intersect $C^5$ at most one
 point
outside of $P$ and thus  we have at most $(C,P') \cong A_1$.
This is a contradiction to $\rho(P',4)>0$.
\qed

\begin{Corollary}\label{NonSimple-list} Assume that $C$ is a reduced
 sextic of torus type with (at least) a non-simple singularity.

\nin
{\rm I.} Assume that $C$ is irreducible. Then
\begin{enumerate}
\item $\De(t)=(t^2-t+1)^2$
if  $\Si(C)$ is either 
$ [C_{3,9},3A_2]$ or $ [B_{3,6},4A_2]$.
\item $\De(t)=(t^2-t+1)$ for other irreducible sextics.
\end{enumerate}
{\rm II.} Assume that $C$ is  not irreducible. Then we have the following 
possibilities.
\begin{enumerate}
\item $\tilde \De(t)=(t^2-t+1)^2$ if   $\Si(C)$ is either 
$ [B_{3,6},4A_2,A_1]$ or $[C_{3,9},3A_2,A_1]$.
\item $\tilde \De(t)=(t^2-t+1)^2(t^2+t+1)$ if $\Si(C)=[2B_{3,6}]$.
\item $\tilde\De(t)=(t^2-t+1)^2(t^2+t+1)(t+1)^2$ 
if $\Si(C)=[B_{3,12}]$
\item $\tilde \De(t)=(t^2-t+1)(t^2+t+1)$ if $\Si(C)$ is one of the
      following.
\[
 [C_{6,6},A_5],\,[C_{6,6},A_5,A_1],\,[C_{6,6}, A_5,2A_1],\,
[C_{6,12}],\,[C_{6,12},A_1],\,[B_{4,6},A_5],\,
[D_{4,7},A_5],\,[Sp_2]
\]
\item $\tilde \De(t)=(t^2-t+1)^4 (t^2+t+1)^4 (t+1)^4$ if $\Si(C)=[B_{6,6}]$.
\item
$\tilde\De(t)=t^2-t+1$
for other reduced sextics of torus curve.
\end{enumerate}
\end{Corollary}
The proof of Corollary is immediate from Theorem \ref{Global-rho-bis} and 
the classification tables in \cite{Pho,Oka-Pho2,Reduced}.
We remark that $\rho(2)\ne  0 $ only for
the configuration $[B_{6,6}]$.

In \cite{Oka-Pho1}, it has been observed that  a tame 
sextics $C$ of torus 
type with $\Si(C)=[3A_2,C_{3,9}]$ is exceptional 
among tame irreducible sextics of torus type in the sense that
$\pi_1(\bfP^2-C)$ is not isomorphic to $\bfZ_2*\bfZ_3$
and its Alexander polynomial is given by $(t^2-t+1)^2$.



\begin{thebibliography}{10}

\bibitem{Artal}
E.~Artal.
\newblock Sur les couples des {Zariski}.
\newblock {\em J.\ Algebraic Geometry}, 3:223--247, 1994.

\bibitem{Esnault}
H.~Esnault.
\newblock Fibre de {M}ilnor d'un c\^one sur une courbe plane singuli\`ere.
\newblock {\em Invent. Math.}, 68(3):477--496, 1982.

\bibitem{LibgoberArcata}
A.~Libgober.
\newblock Alexander invariants of plane algebraic curves.
\newblock In {\em Singularities, Part 2 (Arcata, Calif., 1981)}, pages
  135--143. Amer. Math. Soc., Providence, RI, 1983.

\bibitem{Loeser-Vaquie}
F.~Loeser and M.~Vaqui{\'e}.
\newblock Le polyn\^ome d'{A}lexander d'une courbe plane projective.
\newblock {\em Topology}, 29(2):163--173, 1990.

\bibitem{Merle-Teissier}
M.~Merle and B.~Teissier.
\newblock Conditions d'adjonction, d'apr\`es {D}u {V}al.
\newblock In {\em S\'eminaire sur les singularities des surfaces}, volume 777
  of {\em Lecture Notes in Math.}, pages 230--245. Springer-Verlag, 1977.

\bibitem{Reduced}
M.~Oka.
\newblock Geometry of reduced sextics of torus type.
\newblock {\em math.AG/0203034}.

\bibitem{Two}
M.~Oka.
\newblock Two transforms of plane curves and their fundamental groups.
\newblock {\em J. Math. Sci. Univ. Tokyo}, 3:399--443, 1996.

\bibitem{Okabook}
M.~Oka.
\newblock {\em Non-degenerate complete intersection singularity}.
\newblock Hermann, Paris, 1997.

\bibitem{Oka-Pho2}
M.~Oka and D.~Pho.
\newblock Classification of sextics of torus type.
\newblock {\em math.AG/0201035}.

\bibitem{Oka-Pho1}
M.~Oka and D.~Pho.
\newblock Fundamental group of sextic of torus type.
\newblock In A.~Libgober and M.~Tibar, editors, {\em Trends in Singularities},
  pages 151--180. Birkh\"auser, Basel, 2002.

\bibitem{Pho}
D.~T. Pho.
\newblock Classification of singularities on torus curves of type $(2,3)$.
\newblock {\em Kodai Math. J.}, 24(2):259--284, 2001.

\bibitem{Randell}
R.~Randell.
\newblock Milnor fibers and {A}lexander polynomials of plane curves.
\newblock In {\em Singularities, Part 2 (Arcata, Calif., 1981)}, pages
  415--419. Amer. Math. Soc., Providence, RI, 1983.

\bibitem{Tokunaga-torus}
H.-o. Tokunaga.
\newblock (2,3) torus sextics and the {Albanese} images of 6-fold cyclic
  multiple planes.
\newblock {\em Kodai Math. J.}, 22(2):222--242, 1999.

\bibitem{Tokunaga}
H.-o. Tokunaga.
\newblock Galois covers for {${\frak S}_4$} and {${\frak A}_4$} and their
  applications.
\newblock {\em preprint}, 2000.

\bibitem{Var-zeta}
A.~N. Varchenko.
\newblock Zeta-function of monodromy and {Newton's} diagram.
\newblock {\em Invent. Math.}, 37:253--262, 1976.

\bibitem{Yang}
J.-G. Yang.
\newblock Sextic curves with simple singularities.
\newblock {\em Tohoku Math. J.}, 48(2):203--227, 1996.

\end{thebibliography}
\def\cprime{$'$} \def\cprime{$'$}

\end{document}